\documentstyle{article}

\font\tenmsb=msbm10 scaled\magstep1
\font\sevenmsb=msbm7 scaled\magstep1
\font\fivemsb=msbm5 scaled\magstep1
\newfam\msbfam
\textfont\msbfam=\tenmsb
\scriptfont\msbfam=\sevenmsb
\scriptscriptfont\msbfam=\fivemsb
\def\Bbb#1{{\fam\msbfam\relax#1}}

\def\binom#1#2{{#1\choose#2}}
\def\tbinom#1#2{{\textstyle {#1\choose#2}}}

\def\l{\ell}
\def\C{{\cal C}}

\def\E{{\cal E}}
\def\res{\,\,\cdot\,\,}

\def\I{{{{\scriptstyle I}}}}
\def\NI{{{{\scriptstyle NI}}}}

\def\supp{\mbox{\rm supp}}
\def\Lip{\mbox{Lip}}
\def\iint{\int\!\!\int}

\def\a{\alpha}
\def\be{\beta}
\def\de{\delta}
\def\D{\Delta}
\def\U{\Upsilon}
\def\u{\upsilon}
\def\L{\Lambda}
\def\la{\lambda}
\def\eps{\varepsilon}
\def\g{\gamma}

\def\o{\omega}

\def\z{\zeta}
\def\s{\sigma}
\def\p{\varphi}
\def\x{\xi}

\def\rho{ d}

\def\zbar{\overline{\zeta}}
\def\wbar{\overline{w}}

\def\xibar{\overline{\xi}}

\def\rr{\Bbb R}
\def\nn{\Bbb N}
\def\cc{\Bbb C}

\def\Bbar{\overline{B^n}}

\def\ee{E\times E}

\def\unmig{\frac12}

\def\ad{\be-\frac{n-\la}{p}}

\def\Hgu{H_{\g,1}}
\def\Hgd{H_{\g,2}}
\def\DELTA{\Delta_\g(\z,\x)}
\def\rI{\tau_{{\scriptstyle I}}}
\def\rNI{\tau_{{NI}}}

\def\dbar{\overline{\partial}}
\def\X{\Bbb X}

\def\der#1{\frac{\partial}{\partial{#1}}}
\def\Dt{\tilde{D}}

\def\Bpa{B_\a^p(\mu)}
\def\Bpbe{B_\be^p(\rr^n)}

\def\Lpmu{L^p(d\mu)}
\def\normF{\| f\|_\a}
\def\Lip{\mbox{\rm Lip}}

\def\T{T^{\a}}
\def\Tx{(T_x^\a f)(y)}

\def\Tz{T_\z^\a F(z)}
\def\jet{\{f_j\}_{\o(j)\le\a}}

\def\dmuts{\frac{d\mu(t)\,d\mu(s)}{\mu[t,s]^2}}

\newtheorem{theorem}{Theorem}

\newtheorem{lemma}[theorem]{Lemma}
\newtheorem{proposition}[theorem]{Proposition}
\newtheorem{assumption}[theorem]{Assumption}

\newtheorem{definition}[theorem]{Definition}

\newcommand{\demo}[1]{\par\smallskip \noindent {\bf {#1}:}}
\newcommand{\qed}{\ \hfill\mbox{$\clubsuit$}}

\def\begeq{\begin{equation}}

\title{Extension theorems of Whitney type by the use of integral operators}

\author{Jaume Gudayol\thanks{Partially supported by MEC grant
PB95-0956-c02-01
and CIRIT grant GRQ94-2014. }}

\date{April 16, 1998}

\begin{document}

\maketitle

\begin{abstract}
Given a compact of $\rr^n$, there is always a doubling measure having it as
its support. We use this fact to construct
an integral operator that
extends differentiable functions defined on any compact of $\rr^n$ to the whole
of $\rr^n$. This allows us both to give a new proof of Whitney's extension
theorem and to extend it to Besov spaces defined on arbritrary compact sets
of $\rr^n$. We also modify this operator to obtain, in certain cases, 
holomorphic extensions.
\end{abstract}

\section{Introduction}
Whitney's extension theorem has become a classical tool of analysis, and 
Whitney's proof of it has been widely reused. Our purpose here is to give a 
new proof of it, and a proof which we feel can be more easily adapted to other
contexts.

Let us recall Whitney's theorem
 in its simplest form. Let $E\subset \rr^n$ be a compact set.
Let $\a>0$, and $\{f_j,\,j\in \nn^n,\,|j|\le \a\}$ be a collection of
continuous functions defined on
$E$. For $x\in E$ and $y\in \rr^n$, define:
$$
T_xf(y)=\sum_{|j|\le\a}\frac1{j!} f_j(x)(y-x)^j.
$$
That is, $T_xf(y)$ is the Taylor polinomial of the jet
$\{f_j,\, |j|\le \a\}$ at the point $x$ and evaluated at $y$. 
Assume that
$\Tx$ approximates the jet $\{f_j\}_{|j|\le \a}$ in the same way
in which the Taylor polynomial of a function approximates the function.
What Whitney's extension theorem says is that, under this hypothesis, there
is always a fuction $F(f)$, of class $\C^\a$ in the whole of $\rr^n$,
and of class $\C^\infty$ in $\rr^n\setminus E$,
such that, in a natural sense, the jet $\{f_j\}_{|j|\le \a}$ is the
 restriction to $E$ of
$F(f)$ and its derivatives.

Let us sketch Whitney's proof of this theorem. He begins by covering
$\rr^n\setminus E$ by a collection of balls $\{B(y_i,r_i),i\in\nn\}$ having
suitable properties (this part is what is known as Whitney's covering lemma).
To each of these balls he associates a point $x_i\in E$. Then if
$\{\p_i\}_{i\in\nn}$ is a $\C^\infty$ partition of unity associated to the
covering, Whitney's extension operator is written as:
\begin{equation}
F(f)(y)=\sum_{i\in\nn}\p_i(y)T^\a_{x_i}f(y), \label{Whit}
\end{equation}
where $y\in \rr^n\setminus E$.

Many applications have been found to Whitney's extension operator.
We will explain only two of them. The first one uses it in its original form,
while the second one needs a slight modification of the operator.

Let $B=\{z\in\cc^n, \, |z|<1\}$ be the unit ball of $\cc^n$, and let
$S=\partial B$ be its boundary. Let $E\subset S$ be a closed set. Let $\jet$
be a jet defined on $E$ (we will not be precise about what it means, see
\cite{BrunaOrtega-86} for the details). Let $F(f)$ be the extension to
the whole of $\Bbar$ of $\jet$ defined by \ref{Whit}. Then
if we can solve certain $\dbar$ equation related to $F(f)$ we obtain a
holomorphic function having a prescrived behaviour on $E$.

Now let $E\subset\rr^n$ be a compact set. Then, as seen by Volberg and
Konyagin in \cite{V-K}, there is always a doubling measure $\mu$ having $E$
as its support. With this measure, we can define a Besov space $\Bpa$ on $E$.
In many cases (see \cite{JonssonWallin} or \cite{Jonsson} for details) one
can see that, with some loss of regularity, the restriction from the spaces
$\Bpbe$ to these spaces is well defined. Thus it is natural to ask
for an extension operator from these spaces $\Bpa$ to the spaces $\Bpbe$. In
some particular cases Jonsson and Wallin in
\cite{JonssonWallin} and \cite{Jonsson} gave an extension operator modifying
the definition \ref{Whit}. What they did was to substitute the term
$\T_{x_i}f(y)$
in \ref{Whit} by a term of the form:
$$
\frac1{\mu(B(x_i,Cr_i))}\int_{B(x_i,Cr_i)}T_x^\a f(y)\, d\mu(x),
$$
that is, they subsituted the value at a point (which is not well
defined) by a a mean of the values in a ball.

\smallskip

In both of the previous examples, one can clearly see that it would be more
desirable to have an extension operator of the form:
$$
\E(f)(y)=\int_E K(x,y)\, \Tx\,d\mu(x),
$$
with some suitable kernel and measure. Having it would allow us to restrict
our attention to the existence of a holomorphic kernel, thus avoiding the
$\dbar$ step, in the former case. In the latter, it seems even more natural
to consider extension operators defined by integrals, as we are dealing with
integrable, rather than continuous, functions.

As a matter of fact, such an operator had already been considered in some
particular cases, specially for dealing with analytic functions (see
\cite{Nagel-76}, \cite{tesicarme}, and \cite{BrunaOrtega-91}). But in all the
previous works there where two important restrictions: the
set $E$ was always a variety, and they extended the function, but not (or
with serious restrictions) the derivatives of it.
We present here such an operator. For $x,y\in \rr$, let $\rI(x,y)=|x-y|$, 
and for $z,w\in \cc$, let $\rNI(z,w)=(1-z\wbar)$. Let $h_q^\I$ and $h_q^\NI$ be
given by:
\begin{equation}
h_q(x)=\int_E\tau(x,y)^{-q}d\mu(y),                              \label{defhq}
\end{equation}
for some doubling measure $\mu$. Then the operator has the form:
\begin{equation}
\E(f)(x)=\frac1{h_q(x)}\int_E\frac{\T_y f(x)}{\tau(x,y)^q}d\mu(y),
                                                                \label{defext}
\end{equation}
where $\T_y$ is some Taylor polynomial.
We use this operator to give a new proof of Whitney's extension
theorem. We also use it to prove an extension theorem for Besov spaces
defined on general compact sets of $\rr^n$. We also show that, for subsets $E$
of the unit sphere, in some cases it gives analytic extensions for Lipschitz 
and Besov spaces.

\section{Definitions and statement of results}

\subsection*{The upper dimension of a set}
 Let $(X,\rho)$ be a compact metric space, with diam$(X)<+\infty$.
 For
$x\in X$, $R>0$ and $k\ge1$, let $N(x,R,k)$ be the maximum number of points
lying in
$B(x,kR)$ separated by a distance greater or equal than $R$.
As in \cite{V-K}, we will say that $(X,\rho)\in\Upsilon_\g$ if
 there exists
$C(\g)=C(X,\rho,\g)$ so that, for any $x\in X$ and any $0<R<kR\le 1$,
$$
N(x,R,k)\le C(\g)k^\g.\eqno{(\Upsilon_\g)}
$$
\begin{definition}
Let the upper dimension $\Upsilon(X)$ be:
$$
\Upsilon(X)=\inf\{\g,\,(X,\rho)\in\Upsilon_\g\}.
$$
\end{definition}
This dimension was first introduced by Larman, in \cite{L}, under the
name of uniform metric dimension.

\smallskip

We will say that a probability measure $\mu$ lies in $U_\g=U_\g(X,\rho)$
if there exists $C(\g)$ so that for any $x\in X$ and any $0<R<kR\le 1$,
$$
\mu(B(x,kR))\le Ck^\g \mu(B(x,R)).\eqno{(U_\g)}
$$
Note that, by taking $k=1/R$, $(U_\g)$ implies the weaker condition:
$$
\mu(B(x,R))\ge C R^\g.\eqno{(U'_\g)}
$$
Notice that if for some $\g$, $\mu\in U_\g$, then
$\supp\mu =X$. Moreover, in this case $\mu$ is a doubling
 measure, that is, there exists $C>0$ for which:
$$
\mu(B(x,2r))\le C\mu(B(x,r)).
$$
Let ${\cal U}=\cup_\g U_\g$. It is easily seen (see \cite{V-K}) that ${\cal
U}$ is precisely the set of all doubling measures with support on $X$.

\subsection*{The lower dimension of a set}
Again as in \cite{L},
we will say that $(X,\rho)\in\Lambda_\g$ if there exists
$C(\g)=C(X,\rho,\g)$ so that, for any $x\in X$ and any $0<R<kR\le 1$,
$$
N(x,R,k)\ge C(\g)k^\g.\eqno{(\Lambda_\g)}
$$
Notice that we always have $(X,\rho)\in\Lambda_0$.

\begin{definition}
We define the lower dimension
$\Lambda(X)$ as:
$$
\Lambda(X)=\sup\{\g,\,(X,\rho)\in\Lambda_\g\}.
$$
\end{definition}
This dimension was first defined by Larman
  (\cite{L}) under the name of minimal dimension.

We will say that a doubling measure $\mu$ belongs to
$L_\g=L_\g(X,\rho)$ if
there exists $C(\g)$ so that, for any $x\in X$ and any $0<R<kR\le 1$,
$$
\mu(B(x,kR))\ge Ck^\g \mu(B(x,R)).\eqno{(L_\g)}
$$
As before, by taking $k=1/R$, condition $(\Lambda_\g)$ implies
$$
\mu(B(x,R))\le C R^\g.\eqno{(L'_\g)}
$$
Note that $L_0$ poses no restriction on $\mu\in {\cal U}$.

The following improvement of Volberg and Konyagin's theorem 1 in 
\cite{V-K} can be
found in \cite{Perijo}:
\begin{theorem} \label{mutheorem}
Let $(X,\rho)\in\Upsilon_{\u}\cap\Lambda_{\la}$, for some
$0<\la\le \u<+\infty$.
Then for any $\u'>\u$ and $\la'<\la$ (or $\la'=0$ if $\L(E)=0$) there exists
$\mu\in U_{\u'}\cap L_{\la'}$.
\end{theorem}

\subsection*{The extension operator for the spaces $\Lip_\a(E)$}
Let $E\subset\rr^n$ be a compact set. In $\rr^n$ we consider the metric given 
by $\rho(x,y)=|x-y|$. For a $j\in \nn^n$, let 
$\o(j)=|j|$ be its lenght.
We can define on $E$ the jets $f=\{f_j\}_{\o(j)\le\a}$  to be
collections of continuous functions on $E$, and the Taylor
polynomial of order $\a$ of a jet $f=\jet$ as
\begin{equation}
T^\a_y f(x)=\sum_{\o(j)\le\a}\frac1{j!}f_j(y)w(x,y)^j, \label{Taylor}
\end{equation}
where $w(x,y)=(x-y)$.

Let $\D_{j}(y,x)$ be given by
$$
\D_j(y,x)=f_j(x)-D_x^jT^\a_y f(x). 
$$
We can define then the spaces $\Lip_\a(E)$ as the sets of jets
for which the norm
\begin{equation}
\|f \|_{\Lip_\a(E)}=\sum_{\o(j)\le\a}\left(\|f_j\|_{\infty,E}+
\sup_{x,y\in E}\frac{|\D_{j}(y,x)|}
{\rho(x,y)^{\a-|j|}}\right) \label{Lipschitz}
\end{equation}
is finite.

Between some of this spaces there is a clearly defined derivation
 operator. For each multiindex $j$ we can define:
\begin{eqnarray*}
\Dt^j:\Lip_\a(E) &\mapsto&\Lip_{\a-|j|}(E) \\
\left(f_k\right)_{|k|\le\a}&\to& \left(f_{k+j}\right)_{|k|\le\a-|j|}
\end{eqnarray*}
Notice that 
\begin{equation}
D^j_x T^\a_t f(x)=T^{\a-|j|}_t (\Dt^j f)(x), \label{truc}
\end{equation}
hence it
is natural to define $\Dt^j f=0$ if $|j|>\a$.

With these definitions, we clearly have that
$(\Lip_\a(\rr^n))_{|_E}\subset\Lip_\a(E)$, provided that the
restriction is understood as the induced jet $\{(D^j f)_{|E}\}_{|j|\le \a}$.

Let $E\subset\rr^n$ be a compact set, and let $\u\ge U(E)$. 
Let $\mu\in U_\u (E)$. For such a $\mu$, $q>\u$, 
and $\tau(x,y)=|x-y|$, we define $h_q(x)$ as in \ref{defhq}, and then 
$\E(f)=\E_{q,\a}(f)$ as in \ref{defext}.

Trivially $\E(f)\in {\cal C}^\infty(\rr^n\setminus E)$.
We will see that  $\E(f)\in\Lip_\a (B(0,R))$, for any $R>0$.
We cannot say that
$\E(f)\in\Lip_\a (\rr^n)$, as $\E(f)$ is not necessarily bounded, but this
can be easily obtained by multiplying $\E(f)$ by a suitable support
function. Namely we will prove the following:
\begin{theorem} \label{elprimer}
Let $f\in \Lip_\a(E)$, and $q>\u+\a$. 
Then for any $R>0$  the function $g$ defined by
$$
g(x)=\left\{\begin{array}{ll}
f(x),&\mbox{ if $x\in E$}\\
\E(f)(x),&\mbox{ if $x\notin E$}
\end{array}\right.
$$
lies in $\Lip_\a(B(0,R))$.
\end{theorem}

Let $E$ be any compact subset of $\rr^n$. Some simple considerations that can
be found in \cite{lamevatesi} show that,
for $q>\u>n$, a constant can be found depending only on $R$.

\subsection*{The extension operator for Besov spaces in $\rr^n$}
Let $E\subset \rr^n$ be a compact set and let $\mu\in U_\u(E)\cap L_\la(E)$
be a measure on $E$. A
jet $f=\jet$ will be a collection of functions $f_j\in \Lpmu$.
For such a jet, and $x\in \rr^n$, we
can define $\mu$-a.e. the Taylor polynomial as in \ref{Taylor},
and with it we can define the norm
\begin{equation}
\|f\|_{\Bpa}=\sum_{\o(j)\le\a}\|f_j\|_{\Lpmu}+
\iint_{\ee}\frac{|\D_{j}(t,s)|^p}
{|s-t|^{p(\a-\o(j))-\la}}\dmuts ,           \label{Besov}
\end{equation}
where $\mu[x,y]=\mu(B(x,d(x,y)))$.
Then the Besov space $\Bpa$ will be the set of jets for which this norm is
finite.

If $E=\rr^n$, and $m$ is the Lebesgue measure on $\rr^n$, then
$m(B(t,|t-s|))\approx |t-s|^n$, so that in this case the norm is the
usual for a Besov space, that is
$$
\|f\|_{B^p_\a(\rr^n)}=\sum_{|j|\le\a}\|D^j f\|_{L^p(m)}+
\iint_{\rr^{2n}}\frac{|T^{\a-|j|}_t D^j f(s)-D^j f(t)|^p}
{|s-t|^{p(\a-|j|)+n}}dm(t) dm(s).
$$
Moreover, as seen in \cite{JonssonWallin}, for $f\in\Bpbe$, and
$\a=\ad$, the restriction of $f$ to $\Bpa$ is defined by:
\begin{equation}
D^a f(\x)=\lim_{\de\to 0}\frac1{m(B(\x,\de))}
\int_{B(\x,\de)}D^a f(x)\,dm(x),\label{defrestric}
\end{equation}
for $|a|<\a$, and for $\mu$ almost every $\x\in E$.

In \cite{JonssonWallin} there are a restriction and an extension theorem
for these spaces when there is an $\u$ for which $U_\u\cap L_\u\ne
 \emptyset$, though the restriction theorem
works, in a certain sense, for any doubling measure. In \cite{Jonsson} there
are also restriction and extension theorems, both valid when $\a<1$ (that is,
not involving derivatives).

We are going to prove the following:
\begin{theorem} \label{merde}
Let $E\subset\rr^n$ be a compact set, with $\Upsilon(E)<n$, and $R$ so that
$E\subset B(0,\frac12 R)$. Let $\phi\in{\cal C}^\infty(\rr^n)$
 be a support function for $B(0,2R)$ which is $1$ on $B(0,R)$.
Let $\mu\in U_\u(E)\cap L_\la(E)$, with $\u<n$, and let $\E(f)$ be defined by
 \ref{defext}, for any $q$ large enough.
Assume that $\a\notin \nn$, and let $\be=\a+\frac{n-\la}p$. 
Let $f=\jet\in\Bpa$.
Then $\phi \E(f)\in B^p_\be (\rr^n)$, and for $|a|<\a$, we have that,
in the sense given by \ref{defrestric}, $D^a\E(f)_{|E}=f_a$.
\end{theorem}

\demo{Remark}
The Besov spaces we have just defined are not the standard if $\be\in\nn$, in
the sense that they are always defined by means of first differences, while in
general these spaces are, for some values of the parameters, defined by
means of second differences. Had we used the more usual definition, we would
have found some restrictions in the indices, analogous to those one can find
in extensions from Besov spaces defined on $\rr^n$ to Besov spaces defined on
$\rr^{n+k}$. See theorem \ref{VRB} to see how these restrictions would look
like.

\subsection*{Holomorphic Lipschitz spaces}

Let $B^n$ be the unit ball in $\cc^n$, and let $S$ be the unit sphere.
On $S$ we consider the pseudo-metric $\rho(z,w)=|1-z\wbar|$. 
 As $\rho^{\unmig}$
is a metric, all the previous results concerning the dimensions apply to it
also.
In fact,  we will consider $\rho$ for $z,w\in\Bbar$, 
where, even though it is not a metric,
it satisfies the triangle inequality with a constant.
 We define the holomorphic Lipschitz spaces 
$A_\a(B^n)={\cal H}ol(B^n)\cap {\cal C}^\a(\Bbar)$.
 It is a well known fact that
these spaces are the same as those obtained considering, on $\Bbar$, the metric
$d(x,y)=|x-y|$.

Let $X=\sum_j a_j(z)\der{z_j}+\overline{a_j(z)}\der{\overline{z}_j}$,
where $a_j\in{\cal C}^\infty(\Bbar)$, be a 
vector field. We define its weight $\o(X)$ as $1/2$ if $X$ is 
complex-tangential, i. e. $\sum a_j\overline{z_j} =0$, and $1$ otherwise;
 for a differential operator
$\X=X_1\cdots X_p$ define $\o(\X)=\sum \o(X_j)$.

Let $\z\in S$ be fixed. Let $w_n(z,\z)$ be the normal coordinate, and
let $T^{\cc}_\z(S)$ be coordinated by 
$w_1(z,\z),\dots,w_{n-1}(z,\z)$. For $j\in \nn^n$,
its weight will be $\o(j)=j_n+\unmig (j_1+\dots+j_{n-1})$. With this
weight, and
using the coordinates $w$, we can define the Taylor polinomial of a jet $\jet$
as in \ref{Taylor}. Notice that this polynomial is twice as long in the 
complex-tangential directions (such non isotropic polynomials were
first defined in \cite{FS}). With this polynomial, and if $D^j$ denotes the 
$j$th derivative with respect to the local coordinates $w_1,\dots,w_n$, we
can define the non-isotropic Lipschitz spaces on a closed set $E\subset S$ 
as the sets of continuous jets $\jet$ for which the norm given by 
\ref{Lipschitz} is finite.

Let $E$ be a closed subset of $S$, $\u\ge \Upsilon(E)$, $\la\le \Lambda(E)$ 
and $\mu\in U_\u(E)\cap L_\la(E)$. Let $h_q(z)$ be defined by \ref{defhq},
where $\tau(z,w)=(1-z\wbar)$. For a non-isotropic jet $f=\jet$ whose components
are integrable with respect to $\mu$, let $\E(f)=\E_{q,\a}(f)$ be defined by
\ref{defext}. This extension is well defined and holomorphic wherever 
$h_q(z)\ne 0$. We will assume the following:
\begin{assumption}
\label{cotainfh}
Hereinafter we will assume that $h_q(z)$ also satisfies the bound:
$$
|h_q(z)|\ge C d(z,E)^{-q}\mu(B_z).
$$
\end{assumption}

\demo{Remarks}
It can be checked easily that this is true if $\u< q<1$.
If $E$ is a subset of a complex-tangential variety, then this is true
at least 
for $q<n+\frac14$. This is so because Nagel, in \cite{Nagel-76},
gives the necessary lower bounds for the kernel against which we are
 integrating. 
If $E$ is a complex tangential curve, this is true for all $q$, as 
can be deduced from the results in \cite{tesicarme}.

Under this assumption, we can prove the following:
\begin{theorem}
\label{lemabo}
Let $f\in A_\a(E)$, with $2\a\notin\nn$, and let $\X$ be a differential
operator with  weight $\o(\X)<\a$. Then, for $q>\a+\u$, and whenever 
assumption \ref{cotainfh} holds,
for $z\in\Bbar\setminus E$, 
$$
|\X\E(f)(z)-\E(\X(f)(\res))(z)|\le C(\X)\normF d(z,E)^{\a-\o(\X)}.
$$
Furthermore, if $\o(\X)>\a$ we have
$$
|\X\E(f)(z)|\le C(\X)\normF d(z,E)^{\a-\o(\X)},
$$
and in particular $\E(f)\in A_\a(B)$.
\end{theorem}

\subsection*{Non isotropic Besov spaces}
For $E$ and $\mu$ as in the previous subsection, we can define 
non isotropic Besov spaces $\Bpa$ with respect to $\mu$ as the sets of 
non isotropic jets $\jet$ in $\Lpmu$ for which the norm defined by 
\ref{Besov} is finite. 

If $f\in{\cal C}^1(B^n)$, we define its radial derivative as 
$N f(z)=\sum z_j\der{z_j}f(z)$, and also the derivative $R^1=I+N$.
 Let 
$0<p<+\infty$, $0<q<+\infty$, and $\be\ge0$. Then the Triebel-Lizorkin
space $HF^{p,q}_\be(B^n)$ is the set of holomorphic functions $f$ on $B^n$
so that 
$$
\|f\|_{p,q,\be}^p=\int_S\left(\int_0^1 (1-t^2)^{([\be]+1-\be)q-1}
|R^{[\be]+1}f(tz)|^q dt\right)^{\frac{p}{q}}d\s(z)<+\infty.
$$
For these spaces we prove the following:
\begin{theorem} \label{VRB}
Let $E\subset S$ be a closed set, with $\Upsilon(E)<n$.
Assume that between $\a$ and $\a+(\U(E)-\L(E))/p$
 lies no integer multiple of
$\unmig$.
 Let $\u$, with $\Upsilon(E)\le \u<n$, and $\la \le \Lambda(E)$ 
be close enough so
that between $\a$ and $\a+\frac{\u-\la}p$ lies no integer multiple of $\unmig$.
Let $\mu\in U_\u (E)\cap L_\la(E)$.
Then for $\beta=\a+\frac{n-\la}p$ and
$\E(f)$ defined as in \ref{defext}, with $\jet\in\Bpa$ and
 $q$ large enough, $\E(f)\in HF^{p,1}_\be$, whenever $h_q$ satisfies
assumption \ref{cotainfh}. Also, for $\o(\g)<\a$, we have that,
in the sense given by \ref{defrestric}, $D^a\E(f)_{|E}=f_a$.
\end{theorem}

\section{Technical lemmas}
We begin by
seeing that $h_q$ behaves in a somewhat nice way. Namely, it satisfies
the following:
\begin{proposition} \label{cotesh}
Let $x\notin E$, and let $x_0\in E$ be such that $\rho(x,E)=\rho(x,x_0)$.
Write $B_x=B(x_0,3\rho(x,E))$. Then, for $h_q=h_q^\I,h_q^\NI$ defined by 
\ref{defhq},
\begin{enumerate} 
\item[{\bf (a)}] For $a\ge 0$, $t\in E$, and $q>\u+\a$, 
$$
\int_E\frac{\rho(y,t)^a}{\rho(y,x)^q}d\mu(y)
\le C \rho(t,x)^a\rho(x,E)^{-q}\mu(B_x);
$$
\item[{\bf (b)}] For any $R>0$ and any differential operator $\X$ 
there is a $C=C(\X,R)$ so that if $|x|\le R$,
$$
|\X h_q(x)|\le C \rho(x,E)^{-q-\o(\X)}\mu(B_x);
$$
\item[{\bf (c)}] Under assumption \ref{cotainfh} for $h_q^\NI$,
for any $R>0$ and any differential operator $\X$ 
there is a $C=C(\X,R)$ so that if $|x|\le R$,
$$
|\X \frac1{h_q(x)}|\le C \rho(x,E)^{q-\o(\X)}\mu(B_x)^{-1}.
$$
\end{enumerate}
\end{proposition}

\demo{Proof} We will only prove the bounds for $h_q^\I$, as the bounds for
$h_q^\NI$ are done in exactly the same way, using assumption \ref{cotainfh}
when necessary.
 
For the inequality in {\bf (a)}, the case $a>0$ can be reduced to the 
case $a=0$ applying the triangle inequality to $\rho(y,t)$. In this case, 
we split $h_q(x)$ into the integrals
$$
h_q(x)=\int_{B_x}\frac1{\rho(x,y)^q}d\mu(y)+\int_{E\setminus B_x}
\frac1{\rho(x,y)^q}d\mu(y).
$$
Then the first integral trivially satisfies the upper bound.
As for the second integral, we decompose it into a sum of integrals
over the sets $\{3^j d(x,E)\le |y-x|\le 3^{j+1}d(x,E)\}$, for $j\ge 1$.
On each of these sets, we have that $|y-x|\ge C 3^j d(x,E)$. using it,
and bounding the measure of the set by means of $U_s$, we obtain: 
$$
\int_{E\setminus B_x}\frac1{|x-y|^q}d\mu(y)
\le
 \sum_{j=0}^\infty 3^{-(q-\u)j}d(x,E)^{-q}\mu(B_x)
$$
and, as $q>\u$, this last sum is convergent.

To prove {\bf (b)}, we use that if $y$ lies in a compact set $E$ and
$x\in B(0,R)\setminus E$, then for any $a\in\rr$,
\begin{equation}
|\X(\rho(x,y)^a)|\le C(\X,R,E,a) \rho(x,y)^{a-\o(\X)},\label{VBB}
\end{equation}
so that, using {\bf (a)} we get
$$
|\X h_q(x)|\le C(\X,R) d(x,E)^{-q-\o(\X)}\mu(B_x).
$$

To prove {\bf (c)}, we consider first the case $\o(\X)=0$. Then
it is enough to restrict
the integral to $B_x$, and then use that in this case 
$|x-y|\le 4|x-x_0|=4 d(x,E)$.
Assume now that {\bf (c)} is true whenever $\o(\X)< k$.
Then the use of these bounds in the formula
$$
0=\X (h_q(x)\frac1{h_q(x)})=
\sum_{\o(\X_1)+\o(\X_2)=\o(\X)}
\X_1\left(\frac1{h_q(z)}\right)\X_2\left(h_q(z)\right)
$$
gives us directly the bound for $\o(j)=k$.

\begin{proposition} \label{taylorllest}
Let $x,y\in \rr^n$, and $t\in E$. Then
$$
\T_t f(x)=\sum_{|\l|\le\a}\frac1{\l !}   (x-y)^{\l}
T^{\a-|\l|}_t (\Dt^{\l} f)(y).
$$
\end{proposition}

\demo{Proof}
We decompose $(x-t)$ as $(x-y)+(y-t)$. Thus, if we expand $(x-t)^j$, 
we see that:
$$
\T_t f(x)
=
\sum_{|j|\le\a}\sum_{\l\le j}\frac1{\l !(j-\l)!} f_{\l+(j-\l)}(t)
(x-y)^{\l}(y-t)^{j-\l}.
$$
By rearranging the indices, we get that $|\l|\le\a$, and by writing
$m=j-\l$, we have that $m\ge0$ and $|m|\le \a-|\l|$, so
$$
\T_t f(x)=\sum_{|\l|\le\a}\frac1{\l !}(x-y)^{\l}
\sum_{|m|\le \a-|\l|}\frac1{m!}  f_{\l+m}(t)
(y-t)^{m},
$$
and the inner sum is precisely $T^{\a-|\l|}_t (\Dt^{\l} f)(y)$.\qed

\begin{proposition} \label{conya}
Let $t,s\in E$, and $a,b,c>0$. Let
$$
B_1(t,s)=\{x\in B(0,R),\, \rho(x,s)\le \rho(x,t)\}.
$$
Assume $c-a-b+n<0$. Then if $c-b+n>0$,
$$
\int_{B_1}\frac{\rho(x,E)^c}{\rho(x,t)^a\rho(x,s)^b}dm(x)\le C 
\rho(t,s)^{c-a-b+n}.
$$
\end{proposition}

\demo{Proof}
We split $B_1$
into $A\cup D$, where $A=B_1\cap B(s,\unmig |t-s|)$ and $D=B_1\setminus A$.
Then on $A$, $|x-t|\ge\unmig |t-s|$, and $\rho(x,E)\le \rho(x,s)$, so that 
the integral over $A$ is bounded by
$$
C|t-s|^{-a} \int_{B(s,\unmig |t-s|)} \frac{1}{|x-s|^{b-c}} dm(x)\le
 C|t-s|^{-a}|t-s|^{c-b+n}
$$
whenever $c-b+n>0$.

On the other hand, on $B_1$ we have that $|x-t|^{-1}\le |x-s|^{-1}$, hence the
integral over $D$ is bounded by
$$
  \int_{\rr^n\setminus B(s,\unmig |t-s|)}
 \frac{1}{|x-s|^{a+b-c}} dm(x)\le C |t-s|^{c-a-b+n}
$$
whenever $c-a-b+n<0$. \qed

\begin{lemma} \label{lemashiti}
Let   $0<b<a$. There is a constant C so that for any $z\in \cc$ with
$|z|<1$,
$$
\int_0^1 (1-t)^{b-1}\frac1{|1-tz|^a}dt\le C
\frac1{|1-z|^{a-b}},
$$
and if $b>a$ the integral is bounded by a constant depending only
on $a$ and $b$.
\end{lemma}

\demo{Proof}
The last statement is trivial, we will only prove the first one.
Observe that the integral is always finite between $0$ and $1/2$.
Thus it is enough to bound the integral between $1/2$ and $1$. For
this integral, we use that $|1-tz|\approx|1-z|+1-t$. Then the
computation
of the resulting integral using the change of variables $1-t=s |1-z|$ 
gives us the result.

\section{Proof of theorem 4}
To prove the theorem,
we have to bound $|\T_y g(x)-g(y)|$. To do so, he have to consider
four cases: $x\in E$, and $y\in E$; $x\notin E$, and $y\in E$; $x\in
E$, and $y\notin E$; and $x,y\notin E$. In the first case the bound
comes from the definition, whereas the other possibilities are
considered in the following lemma:
\begin{lemma} \label{THELEMMA}
\begin{enumerate}
\item
Let $x\in B(0,R)$, $y\in E$, $q>\a+\u$ and $|a|\le \a$. Then,
$$
|D^a\E(f)(x)-D^a_x \T_y f(x)|\le C(R) d(x,y)^{\a-|a|}\| f\|_{\a}.
$$
\item
Let $x\in B(0,R)\setminus E$, $y\in E$, $q>\a+\u$ and $|a|\le \a$. Then:
$$
|D^a_y \T_x(\E(f))(y)- f_a(y)|\le C(R) d(x,y)^{\a-|a|}\| f\|_{\a}.
$$
\item
Let $x,y\in B(0,R)\setminus E$, $q>\a+\u$ and $|a|\le\a$. Then
$$
|D^a\T_x \E(f)(y)-D^a \E(f)(y)|\le C(R) \|f\|_{\a} d(x,y)^{\a-|a|}.
$$
\end{enumerate}
\end{lemma}

To prove this lemma, we need the following:
\begin{proposition} \label{smartchange}
Let $x\in \rr^n\setminus E$, $a\in \nn$,
and $q>\u+\a$. then for $x\in B(0,R)$,
$$
|D^a\E_{\a}(f)(x)-\E_{\a-|a|}(\Dt^a f)(x)|\le C(R)
d(x,E)^{\a-|a|}\|f\|_{\Lip_\a(E)},
$$
where, if $|a|>\a$, then $\Dt^a f=0$.
\end{proposition}

\demo{Proof}
We use that 
\begin{equation}
D^a\E_{\a}(f)(x)=
\sum_{k\le a}\tbinom{a}{k}\int_E
 D^{a-k}(\frac1{h_q(x)}\frac{1}{\rho(x,t)^q})D^{k}(\T_t f(x))\, d\mu(t)
\label{llit}
\end{equation}
In this sum, the term corresponding to $k=a$ is $\E_{\a-|a|}(\Dt^a f)(x)$ when
 $|a|\le\a$, whereas if $|a|>\a$ this term does not appear.
Hence to prove the lemma it is enough to bound each of the terms
with $k<a$. We use \ref{truc} and apply to each of these terms proposition
\ref{taylorllest} at a point $x_0\in E$, with $|x-x_0|=d(x,E)$, and it 
decomposes into sums of terms like:
\begin{equation}
(x-x_0)^\l
\int_E D^{a-k}(\frac1{h_q(x)}\frac{1}{|x-t|^q})T_t^{\a-|k|-|\l|}
(\Dt^{k+\l}f)(x_0)\, d\mu(t), \label{ostia}
\end{equation}
with $|\l|\le \a-|k|$.
For each of these terms we are going to use that, as $k<a$,
$0=D^{a-k}(\E(1))$. Then we can substract to it $(x-x_0)^\l f_{k+\l}(x_0)$ 
times this derivative, and \ref{ostia} is the same as 
\begin{equation}
(x-x_0)^\l
\int_E D^{a-k}(\frac1{h_q(x)|x-t|^q})
\D_{k+\l}(t,x_0)\, d\mu(t).  \label{hostia}
\end{equation}
On the other hand, because of proposition \ref{cotesh} and \ref{VBB},
if $\X$ is any differential operator,
\begin{equation}
|\X (\frac1{h_q(x)}\frac{1}{|x-t|^q})|\le  C(R,\X)
d(x,E)^{q-\o(\X)}\mu(B_x)^{-1}\frac{1}{|x-t|^{q}}.        \label{cotanucli}
\end{equation}
Putting this inside the integral in \ref{hostia} (with $\X=D^{a-k}$,
so that $\o(\X)=|a|-|k|$), and using that
$\Dt^{k+\l}f\in\Lip_{\a-|k|-|\l|}(E)$, we get that \ref{hostia} is bounded
by
$$
 C \|f\|_{\a} d(x,E)^{q-|a|+|k|+|\l|}\mu(B_x)^{-1}
\int_E \frac{ |t-x_0|^{\a-|k|-|\l|}}{|x-t|^{q}} d\mu(t),
$$
and again proposition \ref{cotesh} allows us to bound this by 
$C\rho(x,E)^{\a-|a|}\|f\|_\a$.\qed

\demo{Proof of part 1 of lemma \ref{THELEMMA}}
Because of \ref{truc}, and using proposition \ref{smartchange},
$$
|D^a\E(f)(x)-D^a_x \T_y f(x)|\le
C d(x,E)^{\a-|a|}\|f\|_{\a}+|\E(\Dt^a f)(x)-T^{\a-|a|}_y (\Dt^a f)(x)|.
$$
But if $y\in E$, then $d(x,E)\le d(x,y)$ so we only need to bound the second
term of the inequality. Thus we have to bound
\begin{equation}
\frac1{h_q(x)}\int_E
\frac{|T^{\a-|a|}_t (\Dt^a f)(x)-T^{\a-|a|}_y (\Dt^a f)(x)|}{|x-t|^q}d\mu(t).
\label{culmireya}
\end{equation}
Using proposition \ref{taylorllest}, and as 
$\Dt^{a+\l}f\in\Lip_{\a-|a|-|\l|}(E)$, 
we get that
$$
|T^{\a-|a|}_t (\Dt^a f)(x)-T^{\a-|a|}_y (\Dt^a f)(x)|\le
 C \sum_{|\l|\le \a-|a|}|x-y|^{|\l|} |t- y|^{\a-|a|-|\l|}\|f\|_\a.
$$ 
If we put this inside the integral in \ref{culmireya}, we see that
\ref{culmireya} is, because of proposition \ref{cotesh}, less than
$$
C\frac1{h_q(x)}\|f\|_{\a} \sum_{|\l|\le \a-|a|}d(x,y)^{|\l|}\int_E
\frac{|y-t|^{\a-|a|-|\l|}}{|x-t|^q}d\mu(t)\le C \|f\|_{\a}d(x,y)^{\a-|a|},
$$
as we wanted to see.\qed

\demo{Proof of part 2 of lemma \ref{THELEMMA}}
We use \ref{truc} again, and developing the Taylor polynomial that we obtain 
gives us that $|D^a_y \T_x(\E(f))(y)-\Dt^a f(y)|$ is bounded by:
\begin{eqnarray}
 && \sum_{|\l|\le \a-|a|}\frac1{\l !}|x-y|^{|\l|}
|D^{a+\l}\E(f)(y)-\E(\Dt^{a+\l}f)(y)|+ \nonumber
\\
&+& | \sum_{|\l|\le \a-|a|}\frac1{\l !}(x-y)^{\l}
\E(\Dt^{a+\l}f)(y)-\Dt^a f(y)|. \label{pepito}
\end{eqnarray}
Now proposition \ref{smartchange} shows that the former term is bounded by:
$$
 \sum_{|\l|\le \a-|a|}\frac1{\l !}d(x,y)^{|\l|}d(x,E)^{\a-|a|-|\l|}
\|f\|_{\a}\le C d(x,y)^{\a-|a|}\|f\|_{\a}.
$$
To bound the latter term, we expand $\E(\Dt^{a+\l}f)(y)$, enter the sum inside 
the integral and apply  proposition \ref{taylorllest}, and thus
\begin{equation}
\sum_{|\l|\le \a-|a|}\frac1{\l !}(x-y)^{\l}
\E(\Dt^{a+\l}f)(y)
=\frac1{h_q(x)}\int_E\frac1{|x-t|^q}
T_t^{\a-|a|}(\Dt^a f)(y)\,d\mu(t).\label{mandra}
\end{equation}
Therefore, \ref{pepito} is bounded by  
$$
 \frac1{h_q(x)}\int_E\frac{|T_t^{\a-|a|}(\Dt^a f)(y)-\Dt^a f(y)|}{|x-t|^q}
\,d\mu(t).
$$
Now as $\Dt^a f\in\Lip_{\a-|a|}$, and because of part {\bf (a)} in 
proposition \ref{cotesh},  
this is bounded by $ C \|f\|_{\a} d(x,y)^{\a-|a|}$, as we wanted to see.

\demo{Proof of part 3 of lemma \ref{THELEMMA}}
To prove 3, we will split $B(0,R)^2$ into the sets
\begin{equation}
A_1=\{(x,y), \, \rho(x,y)\le \frac14 \max\{\rho(x,E),\rho(y,E)\}\,\},
\label{apa}\end{equation}
and $A_2=B(0,R)^2\setminus A_1$.
Then on $A_1$, for all $\x=a x+(1-a) y$ with $0\le a \le1$ we have
$d(\x,E)\ge 3 d(x,y)$. 
Let $j$ be any multiindex with $|j|=[\a]+1$. As $\E(f)$ 
is ${\cal C}^\infty$ outside $E$, the mean value theorem shows that
\begin{equation}
|D_x^a(\T_y\E(f)(x)-\E(f)(x))| 
\le C\sup_{|j|=[\a]+1}\sup_{\x\in [x,y]}|D^j\E(f)(\x)|\,|x-y|^{[\a]+1-|a|}.
\label{meanvalue}
\end{equation}
Now proposition
\ref{smartchange}, when applied to $D^j\E(f)(\x)$, allows us to bound 
\ref{meanvalue} by $C(R) \|f\|_{\a}d(x,y)^{\a-|a|}$.

We have to bound the same difference for $(x,y)\in A_2$. To do so,
we use \ref{truc},
expand $T_y^{\a-|a|}D^a\E(f)$ and approximate $D^{a+j}\E(f)$ by 
$\Dt^{a+j} f$. Then the difference we have to bound can be 
split into tree terms, namely:
\begin{eqnarray}
&& |\sum_{|j|\le \a-|a|}\frac1{j!}(x-y)^j(D^{a+j}\E(f)(y)-
\E(\Dt^{a+j} f)(y))|+\label{primerterme} 
\\
&& \qquad
+|\sum_{|j|\le \a-|a|}\frac1{j!}(x-y)^j\E(\Dt^{a+j}f)(y)-\E(\Dt^a f)(x) |
+\label{segonterme} 
\\
&&\qquad + |\E(\Dt^a f)(x)-D^a\E(f)(x)|. \label{tercerterme}
\end{eqnarray}
For the term \ref{tercerterme} we use proposition \ref{smartchange} and
bound it by:
$$
|\E(\Dt^a f)(x)-D^a\E(f)(x)|\le C(R) \|f\|_{\a}d(x,E)^{\a-|a|}\le
C(R) \|f\|_{\a}d(x,y)^{\a-|a|}.
$$
In the same way, \ref{primerterme} is bounded by:
\begin{eqnarray*}
&& 
 C(R) \|f\|_{\a}\sum_{|j|\le \a-|a|}\frac1{j!}d(x,y)^{|j|}
d(x,E)^{\a-|a|-|j|}\le C(R) \|f\|_{\a}d(x,y)^{\a-|a|}.
\end{eqnarray*}

To bound \ref{segonterme}, 
we use proposition \ref{taylorllest} in the same way as in
\ref{mandra}  and then multiply by
$(|x-t|^q h_q(x))^{-1}$ and integrate against $d\mu(t)$. Thus we obtain:
$$ 
\sum_{|j|\le \a-|a|}\frac{(x-y)^j}{j!}\E(\Dt^{a+j}f)(y)
=\iint_{\ee} \frac
{\sum_{|\l|\le \a-|a|}\frac1{\l !}(x-s)^{\l} f_{a+\l}(s)}
{h_q(y)h_q(x)|x-t|^q |y-s|^q}d\mu(s)\,d\mu(t).
$$ 
On the other hand, and also because of proposition \ref{taylorllest}, 
and proceeding as before,
$$
\E(\Dt^a f)(x)=
\iint_{\ee}\frac
{\sum_{|\l|\le \a-|a|}\frac1{\l !}(x-s)^{\l}
T_t^{\a-|a|-|\l|}(\Dt^{a+\l} f)(s)}
{h_q(y)h_q(x)|x-t|^q|y-s|^q}d\mu(t)d\mu(s).
$$
If we add up these two facts, we get that \ref{segonterme} can be bounded by:
$$
\iint_{\ee}\frac
{\sum_{|\l|\le \a-|a|}\frac1{\l !}|x-s|^{|\l|}| f_{a+\l}(s)-
T_t^{\a-|a|-|\l|}(\Dt^{a+\l} f)(s)|}
{h_q(y)h_q(x)|x-t|^q |y-s|^q}d\mu(s)\,d\mu(t).
$$
Now as $\Dt^{a+\l} f\in \Lip_{\a-|a|-|\l|}(E)$, this can be bounded by:
$$
 C\|f\|_{\a}\sum_{|\l|\le \a-|a|}\frac1{\l !}
\frac1{h_q(y)h_q(x)}\iint_{\ee} \frac
{|x-s|^{|\l|}|s-t|^{\a-|a|-|\l|}}
{|x-t|^q |y-s|^q}d\mu(s)\,d\mu(t).
$$
If we use now that $|x-s|^{p}\le C_p(|x-t|^p+|t-s|^p)$, and then that 
$|t-s|^{p}\le C_p(|x-t|^p+|x-y|^p+|y-s|^p)$, and apply proposition 
\ref{cotesh} to each of the integrals that result from it,
 we can bound each of these integrals by sums of
terms of the form
$$
d(x,y)^{a_1}d(x,E)^{a_2}d(y,E)^{a_3}
$$
where $a_1,a_2,a_3\ge 0$ and
$a_1+a_2+a_3=\a-|a|$. But as both $d(x,E)$ and $d(y,E)$ are less than
a constant times $d(x,y)$, we are done.\qed

\section{Proof of theorem 5}
\begin{lemma} \label{ojala}
Under the assumptions of theorem \ref{merde}, if $a$ is any multiindex,
$$
\|\E(\Dt^a f)(x)\|_{L^p(B(0,R))}\le C \|f\|_{B^p_\a(\mu)}.
$$
\end{lemma}

\demo{Proof}
If $|a|>\a$, then $\Dt^a f=0$, so we only have to consider $|a|\le \a$. On the
other hand, entering the modulus inside the integral, we see that what we 
have to bound is:
$$
C(R)\sum_{|k|\le \a-|a|}
\int_{B(0,R)}\frac1{h_q(x)^p}\left(\int_E \frac{|\Dt^{a+k} f(t)|}
{|x-t|^q} d\mu(t)\right)^p dm(x).
$$
If we apply H\"older's inequality  and bound the integral not 
containing $|\Dt^{a+k} f(t)|$ using proposition \ref{cotesh} and
$U_\u'$, we can estimate this last integral by 
$$
 C \sum_{|k|\le \a-|a|}
\int_E\int_{B(0,R)} \frac{|f_{a+k} (t)|^p}{|x-t|^\u}dm(x) d\mu(t)
\le C \sum_{|k|\le \a-|a|} \| f_{a+k}\|_{\Lpmu},
$$
as $\u<n$.\qed

\begin{lemma} \label{sandra}
Under the assumptions of theorem \ref{merde},
if $a$ is a multiindex with $|a|\le\be$, then
$$
\| D^a \E(f)(x)-\E(\Dt^a f)(x)\|_{L^p(B(0,R))}\le C\| f\|_{\Bpa}.
$$
\end{lemma}

\demo{Proof}
Let $s\in E$. Proceeding for this $s$ as we procceded with $x_0$ in the
the proof of proposition \ref{smartchange} we see that 
it is enough to bound terms like \ref{hostia}. But using \ref{cotanucli} 
allows us to bound \ref{hostia} by 
\begin{equation}
 C |x-s|^{|\l|}\frac{d(x,E)^{q-|a|+|k|}}{\mu(B_x)}
\int_E \frac{|\D_{k+\l}(t,s)|}{|x-t|^q} d\mu(t),
\label{star}\end{equation}
with $k\le a$ and $|k|+|\l|\le \a$. We will write $m=k+\l$.
We now integrate against $|x-s|^{-q-|\l|}d\mu(s)$ and divide the result by
$h_{q+|\l|}(x)$. If we apply the bounds in \ref{cotesh} to  $h_{q+|\l|}(x)$,
then the previous term is bounded by
$$
 C \frac{d(x,E)^{2q-|a|+|m|}}{\mu(B_x)^2}
\iint_{\ee} \frac{|\D_{m}(t,s)|}{|x-t|^q |x-s|^{q}} d\mu(t)d\mu(s)
=I(x).
$$
Let $0<A,B<q$ be such that $(q-A)p'>\u $ and $(q-B)p'>\u$ (we can always 
choose such an $A$ and $B$, for $q$ large enough).
 By H\"older's inequality, and using {\bf (a)} in proposition \ref{cotesh}  
in the integral not containing $\D_{m}(t,s)$, we have that
\begin{equation}
 I(x)^p \le 
C \frac{d(x,E)^{p(A+B-|a|+|m|)}}{\mu(B_x)^{2}}
\iint_{\ee} \frac{|\D_{m}(t,s)|^p d\mu(t)d\mu(s)}
{|x-t|^{Ap} |x-s|^{Bp}}.
\label{permestard}
\end{equation}
If we now use that $d(x,E)$ is bounded above and that
$|a|\le \be=\a+\frac{n-\la}p$, and apply Fubinni's theorem, we can bound
$\int_{B(0,R)}I(x)^p\,dm(x)$ by
\begin{equation}
\iint_{\ee}|\D_{m}(t,s)|^p\int_{B(0,R)}
\frac{d(x,E)^{p(A+B-\a+|m|)-n+\la}}{|x-t|^{Ap} |x-s|^{Bp}\mu(B_x)^{2}}
dm(x)d\mu(t)d\mu(s).                      \label{isa}
\end{equation}

Fix $(t,s)\in \ee$. Let us split $B(0,R)$ into the sets
\begin{equation}
B_1(t,s)=\{x\in B(0,R),\, |x-s|\le |x-t|\}
\quad\mbox{ and }\quad B_2=B(0,R)\setminus B_1        \label{split}
\end{equation}
It is enough to bound the integral over $B_1$, as that over $B_2$ is bounded
similarly, changing the roles of $s$ and $t$.
On $B_1$, $|t-s|\le 2|x-t|$, so that by $U_\u$ and $L_\la$,
$$
\mu(B(t,|t-x|))\ge C
\mu(B(t,|t-s|\,2\frac{|t-x|}{|t-s|}))\ge C
\frac{|t-x|^\la }{|t-s|^\la }\mu(B(t,|t-s|)).
$$
On the other hand, $B(t,|t-x|))\subset B(x_0, 3|t-x|)$, so that
$$
\mu(B(t,|t-x|))\le C 
\mu(B(x_0,3d(x,E)\frac{|t-x|}{d(x,E)}))\le C
\frac{|t-x|^\u}{d(x,E)^\u}\mu(B_x).
$$
If we add up the two facts, we get that
\begin{equation}
\mu(B_x)\ge C\frac{d(x,E)^\u}{|x-t|^{\u-\la}|t-s|^\la}\mu[t,s].
                                                          \label{cotamuBx}
\end{equation}
Hence, the part of the inner integral in \ref{isa} corresponding to $B_1$ is
bounded by
$$
 C \frac{|t-s|^{2\la}}{\mu[t,s]^2}\int_{B(0,R)}
\frac{d(x,E)^{p(A+B-\a+|m|)-n+\la-2\u}}{|x-t|^{Ap-2(\u-\la)}|x-s|^{Bp}}
dm(x).
$$
We apply proposition \ref{conya} to this integral, something we can do
if we chose $A$ and $B$ properly, as $|m|=|k|+|\l|<\a$.
 Then we get that it is bounded by
$\mu[t,s]^{-2}$ times $|t-s|$ to the power
$-p(\a-|k|-|\l|)+\la<0$, as we wanted to see. \qed

The following lemma insures that $\E(f)$ interpolates $f$:
\begin{lemma} \label{viscaelbarsa}
For $|j|<\a$, and $\mu$-a.e. $\x\in E$,
$$
\lim_{\de\to 0}\frac1{m(B(\x,\de))}
\int_{B(\x,\de)}|D^j \E(f)(x)-f_j(\x)|^p\,dm(x)=0.
$$
\end{lemma}

\demo{Proof}
In the proof of the previous lemma, we saw that
$|D^j\E(f)-\E(\Dt^j f)|$  could be bounded by terms like 
\ref{permestard}. Now if we use that in this case we have that
$|a|\le\a$ instead of $|a|\le \be$, and proceed as there, we get that 
\ref{permestard} can be bounded by $d(x,E)^{\frac{n-\la}p}$ times some
integrable function. Thus this tends to zero if we calculate the mean
value of this on the ball $B(\x,\de)$ and let $\de\to 0$. Hence, what
we have to see is that
$$
\lim_{\de\to 0}\frac1{m(B(\x,\de))}
\int_{B(\x,\de)}|\E(\Dt^j f)(x)-f_j(\x)|^p\,dm(x)=0.
$$
Now  if we apply proposition \ref{taylorllest} to 
$\E(\Dt^j f)(x)-f_j(\x)$ we see
that it decomposes into:
\begin{eqnarray*}
&& \frac1{h_q(x)}\int_E\frac1{|x-y|^q}(T^{\a-|j|}_y(\Dt^j
f)(\x)-f_j(\x)) \,d\mu(y)
+\\
&+& \frac1{h_q(x)}\int_E\frac1{|x-y|^q}\sum_{0<|k|<\a-|j|}\frac1{k!}
(x-\x)^k
T^{\a-|j|-|k|}_y(\Dt^{j+k} f)(\x)\,d\mu(y)
\end{eqnarray*}
The latter term is nothing but a sum of terms like 
$(x-\x)^k \E(\Dt^{k+j} f)(\x)$, with $|k|>0$. Thus this is integrable
w.r.t. $\x\in  E$, thus finite a.e. Then this part is bounded by 
$\de^{|k|}$ times some integrable function, thereby tending to $0$
with $\de$.

The former term is bounded by something like \ref{star}, where $\l=0$,
and $k=a=j$. We raise it to the power $p$, apply H\"older's inequality
to what we obtain, and then use
part {\bf (a)} of proposition \ref{cotesh} to bound the integral not
containing $|\D_j(y,\x)|$. We use also that $d(y,\x)^{-\la}\le
\mu[y,\x]^{-1}$. 
Thus this term is bounded by:
\begin{equation}
\frac{d(x,E)^{\u}d(x,\x)^{(\a-|j|)p+\la}}{\mu(B_x)}
\int_E \frac{|\D_{j}(t,\x)|^p}{d(x,y)^\u d(\x,y)^{(\a-|j|)p}} 
\frac{d\mu(t)}{\mu[\x,y]}.                               \label{integrable}
\end{equation}
To bound this integral, we split $E$ into the sets 
$E_1=\{y,\, d(y,x)\ge 2 d(x,\x)\}$ and $E_2$ its complementary.

To bound the part of \ref{integrable} corresponding to $E_1$, we first
use the equivalent of \ref{cotamuBx}. Next we use that on $E_1$ we
have $d(x,y)^{-\la}\le d(x,\x)^{-\la}$. Hence, and as $d(x,\x)\le \de$,
this part is bounded by $\de^{(\a-|j|)p}$ times some integrable (and
thereby finite $\mu$-a.e.) function of $\x$.

To bound the part of \ref{integrable} corresponding to $E_2$, we use
again the equivalent of \ref{cotamuBx} (which is different that for
$E_1$). This time we bound $d(x,\x)$ by $\de$ first, and next we
integrate over $B(\x,\de)$ and apply Fubini's theorem. But as 
$y\in E_2$, we also have $y\in B(\z,3\de)$ and thus the inner
integral can be bounded by $\de^{n-\u}$. Hence we finaly obtain the
same bound as before, so we are done. \qed

To prove the theorem, we will split $B(0,R)^2=A_1\cup A_2$ as in \ref{apa}.
Then we have:
\begin{lemma} \label{cotallesta}
Let $a$ be a multiindex with $|a|\le\be$. Then
$$
\iint_{A_2}\frac{|D^a\E(f)(x)-\E(\Dt^a f)(x)|^p}{|x-y|^{p(\be-|a|)+n}}
dm(x)dm(y)\le C\|f\|_{\Bpa}
$$
\end{lemma}

\demo{Proof}
In the proof of lemma \ref{sandra} we have seen that the difference
$|D^a\E(f)(x)-\E(\Dt^a f)(x)|^p$ 
is bounded by sums of terms like $I(x)^p$, as defined in \ref{permestard}.
On the other hand
$$
\displaylines{
\iint_{A_2} \frac{I(x)^p}{|x-y|^{p(\be-|a|)+n}}dm(x)dm(y)\le \hfill\cr
\hfill \le C \int_{B(0,R)} I(x)^p
\int_{\{y,\,2R\ge |x-y|\ge \frac13 d(x,E) \}}
\frac{dm(y)}{|x-y|^{p(\be-|a|)+n}} dm(x)\le\hfill\cr
\le C \int_{\ee}|\Delta_{k+\l}(t,s)|^p \int_{B(0,R)^2}
\frac{d(x,E)^{p(A+B-\a+|k|+|\l|)-n+\la-2\u}}
{|x-t|^{Ap} |x-s|^{Bp}\mu(B_x)^{2}}
dm(x)d\mu(t)d\mu(s),
}$$
and these integrals have already been bounded in lemma \ref{sandra}.

\demo{Proof of the theorem}
We have to see that the integral
$$
\iint_{B(0,R)^2}
\frac{|T_x^\be \E(f)(y)-\E(f)(y)|^p}{|x-y|^{\be p+n}}dm(x)dm(y)
$$
is bounded. The $L^p$ norms have already been bounded by the previous
lemmas.

As we have already said, we will split $B(0,R)^2=A_1\cup A_2$ 
as in \ref{apa}.
We will begin by bounding the integral over $A_1$. In $A_1$ we can apply
the mean value theorem as in \ref{meanvalue}. Thus we have to bound
$|D^a\E(f)(\x)|$, with $|a|=[\be]+1$ and
$d(\x,E)\ge 3|x-y|$, and in particular
$\x\notin E$ (recall that, as $\U(E)<n$, the measure of $E\times E$ is
zero). 
On the other hand, as we have seen in the proof of lemma \ref{sandra},
 $|D^a\E(f)(\x)|$ is bounded by
sums of terms like $I(\x)$, and $I(\x)^p$ is bounded in \ref{permestard}.

As we have seen in \ref{cotamuBx}, either
$$
\mu(B_{\x})\ge C\frac{d(\x,E)^\u}{|\x-t|^{\u-\la}|t-s|^\la}\mu[t,s]
\quad \mbox{ or }\quad
\mu(B_{\x})\ge C\frac{d(\x,E)^\u}{|\x-s|^{\u-\la}|t-s|^\la}\mu[t,s].
$$
Recall that $|x-y|\le\frac14 \max \{d(x,E),d(y,E)\}$. Hence
$|x-t|\approx |\x-t|\approx|y-t|$, and the same is true for $s$ instead of $t$.
Also $d(x,E)\approx d(\x,E)\approx d(y,E)$. Therefore, $I(\x)^p$ is bounded by
$$
d(y,E)^{p(A+B-|a|+|k|+|\l|)-2\u}
\iint_{\ee}\frac{|\D_{k+\l}(t,s)|^p|t-s|^{2\la}}
{|x-t|^{Ap-\u+\la}|y-s|^{Bp-\u+\la}}\frac{d\mu(t)d\mu(s)}{\mu[t,s]^2}.
$$
Thus the integral over $A_1$ is bounded by
\begin{equation}
 \iint_{\ee}|\D_{k+\l}(t,s)|^p|t-s|^{2\la}J(s,t)
\frac{d\mu(t)d\mu(s)}{\mu[t,s]^2},                 \label{zira}
\end{equation}
where
$$
J(s,t)=\iint_{A_1}
\frac{|x-y|^{(|a|-\be)p-n}d(y,E)^{p(A+B-|a|+ |k|+|\l|)-2\u}}
{|x-t|^{Ap-\u+\la}|y-s|^{Bp-\u+\la}}dm(x)dm(y).
$$
For $J(s,t)$ we use that $|x-t|\approx |y-t|$. Hence $J(s,t)$ is bounded
by: 
$$ 
\int_{B(0,R)}
\frac{d(y,E)^{p(A+B-|a|+|k|+|\l|)-2\u}}
{|y-t|^{Ap-\u+\la}|y-s|^{Bp-\u+\la}}
\int_{B(y,\frac14 d(y,E))}|x-y|^{(|a|-\be)p-n}dm(x)dm(y).
$$
As $|a|-\be>0$, the inner integral is bounded by $d(y,E)^{(|a|-\be)p}$.
Then if we  apply
proposition \ref{conya} to the integral with respect to $y$,
 we get that $J(s,t)$ is bounded by
$|t-s|^{-p(\a-|k|-|\l|)-\la}$
and this, when used in \ref{zira}, says that \ref{zira} is bounded by 
$ \|f\|_{\Bpa}$,
as we wanted to see.

\smallskip

We have to bound now the integral over $A_2$. But if we approximate
each of the derivatives in $T_x^\be \E(f)$ by $\E(\Dt^j f)(x)$, then
the integral over $A_2$ is bounded by:
 \begin{eqnarray}
&&\sum_{|j\le\be} \iint_{A_2} \frac
{|D^j \E(f)(x)-\E(\Dt^j f)(x)|^p}
{|x-y|^{(\be-|j|) p+n}}dm(x)dm(y)+\nonumber
\\
&+&\iint_{A_2} \frac
{|\sum_{|j|\le\a}\frac1{j!}(y-x)^j\E(\Dt^j f)(x)-\E(f)(y)|^p}
{|x-y|^{\be p+n}}dm(x)dm(y), \label{COOL}
\end{eqnarray}
and we only have to bound the latter term, as the former is already bounded
by lemma \ref{cotallesta}.

Applying proposition \ref{taylorllest} twice, the first to reagrupate
terms, and the seccond with some
$s\in E$, leads to:
$$
\sum_{|j|\le\a}\frac{(y-x)^j}{j!}\E(\Dt^j f)(x)
=\frac1{h_q(x)}\int_E \sum_{|k|\le\a}\frac{(y-s)^k}{k!|x-t|^q}
T_t^{\a-|k|}(\Dt^k f)(s)d\mu(t),
$$
Then, if we integrate with respect so $|y-s|^{-q}d\mu(s)$
and divide by $h_q(y)$, we get that this sum is equal to
$$ 
\sum_{|k|\le\a}\frac1{k!}\frac1{h_q(x)h_q(y)}\iint_{\ee} \frac
{(y-s)^k T_t^{\a-|k|}(\Dt^k f)(s)}{|x-t|^q |y-s|^q}
d\mu(t)d\mu(s).
$$ 
Likewise
$$
\E( f)(y)=
\sum_{|k|\le\a}\frac1{k!}\frac1{h_q(x)h_q(y)}\iint_{\ee}
\frac{(y-s)^k f_k(s)}{|x-t|^q|y-s|^q}
d\mu(s)d\mu(t).
$$
Therefore the difference appearing inside the integral in \ref{COOL}
is bounded by sums of terms of the form 
$$
A_k(x,y)= \frac1{h_q(x)h_q(y)}\iint_{E\times E} \frac
{| T_t^{\a-|k|}(\Dt^k f)(s)-f_k(s)|}{|x-t|^q |y-s|^{q-|k|}}
d\mu(t)d\mu(s)
$$
For $|k|\le\a$. 
Let $k$ be fixed, and let $0<a<q$ and $0<b<q-|k|$ be some numbers to
 be chosen later. Because of H\"older's inequality and part {\bf(a)} 
in proposition \ref{cotesh},
$$
A_k(x,y)^p\le
C \frac{d(x,E)^{ap}d(y,E)^{(b+|k|)p}}{\mu(B_x)\mu(B_y)}
\iint_{\ee}\frac{|\D_{k}(t,s)|^p}{|x-t|^{ap}|y-s|^{bp}}d\mu(t)d\mu(s).
$$
Thus,
\begin{eqnarray*}
&&\iint_{A_2} \frac{A_k(x,y)^p}{|x-y|^{\be p+n}}dm(x)dm(y)\le
C \iint_{\ee}\frac{|\D_{k}(t,s)|^p}{|s-t|^{(\a-|k|)p-\la}}
\frac{d\mu(t)d\mu(s)}{\mu[t,s]^2}\times\\
&&\times\sup_{s,t\in E}
\iint_{A_2} \frac{|s-t|^{(\a-|k|)p-\la}d(x,E)^{ap}d(y,E)^{(b+|k|)p}\mu[t,s]^2}
{|x-y|^{\be p+n}|x-t|^{ap}|y-s|^{bp}\mu(B_x)\mu(B_y)}dm(x)dm(y).
\end{eqnarray*}

To finish the proof, we only have to see that this supremum is finite.
To do so, we use that, on $A_2$,
 $d(y,E)\le 4|x-y|$, so that our integral
is bounded by
$$
\iint_{A_2} \frac{|s-t|^{(\a-|k|)p-\la}d(x,E)^{ap}d(y,E)^{bp}\mu[t,s]^2}
{|x-y|^{(\be -|k|)p+n}|x-t|^{ap}|y-s|^{bp}\mu(B_x)\mu(B_y)}dm(x)dm(y)
$$
and this integral is symetric with respect to swapping $x$ for $y$ and
$s$ for $t$.

Let $B_{1,1}=B(s,\frac13 |t-s|)$ and $B_{1,2}=B_1\setminus B_{1,1}$.
Let $B_{2,1}=B(t,\frac13 |t-s|)$ and $B_{2,2}=B_2\setminus B_{2,1}$.
We split the integral over $B$ into the integrals over the sets
$$
B_{i,j}^{k,\l}=\{(x,y)\in B,\,\, x\in B_{i,j},\,y\in B_{k,\l}\},
$$
where $i,j,k,\l\in \{1,2\}$.
Then we will have to bound 16 integrals. Because of the simetries of the
integrand we only have to estimate the integrals over $B_{1,1}^{1,1}$,
$B_{1,1}^{1,2}$,
$B_{1,1}^{2,1}$, $B_{1,1}^{2,2}$,  $B_{1,2}^{1,2}$,  and $B_{1,2}^{2,2}$.
All of them are done in more or less the same way, using 
\ref{cotamuBx} suitably (the fact that we have to use \ref{cotamuBx} each
time in a different way is what makes necessary to consider so many cases).
We will bound the integral over $B_{1,2}^{1,2}$, which is perhaps the more
delicate, and the others are done similarly.

for the integral over $B_{1,2}^{1,2}$, we have that $x\in B_1$ and
$y\in B_1$, so that
$$
\frac{\mu[t,s]^2}{\mu(B_x)\mu(B_y)}\le C\frac{|x-t|^{\u-\la}|y-t|^{\u-\la}}
{d(x,E)^\u d(y,E)^\u}|t-s|^{2\la}.
$$
Thus what we have to bound is
$$
|s-t|^{(\a-|k|)p+\la}
\iint_{B_{1,2}^{1,2}} \frac{d(x,E)^{ap-\u}d(y,E)^{bp-\u}|y-t|^{\u-\la}}
{|x-y|^{(\be -|k|)p+n}|x-t|^{ap-\u+\la}|y-s|^{bp}}dm(x)dm(y).
$$
We will use only that $|x-y|\ge \frac14 d(x,E)$ and $|x-y|\ge \frac14 d(y,E)$.
As $(\be-|k|)p+n=(\a-|k|)p+2n-\la$, and $(\a-|k|)p>0$ because $\a\notin\nn$,
we can write $(\a-|k|)p=\de_1+\de_2$ with $\de_1,\de_2>0$. Then the integral is
bounded by
\begin{eqnarray*}
&&|s-t|^{(\a-|k|)p+\la}
\int_{B_{1,2}} \frac{d(x,E)^{ap-\u-\de_1-n+\la}}
 {|x-t|^{ap-\u+\la}}dm(x)\times
\\
&&\qquad\qquad\times\int_{B_{1,2}} \frac{d(y,E)^{bp-\u-\de_2-n}|y-t|^{\u-\la}}
{|y-s|^{bp}}dm(y)\le
\\
&&\le C |s-t|^{(\a-|k|)p+\la}\int_{B_{1,2}} \frac{dm(x)}{|x-t|^{n+\de_1}}
 \int_{B_{1,2}}\frac{|y-t|^{\u-\la}}
{|y-s|^{n+\de_2+\u}}dm(y)\le C.
\end{eqnarray*}

\section{Proof of theorem 7}
The following proposition can be found in \cite{BrunaOrtega-86},
 lemma 2.1, and will be useful to us:
\begin{proposition}
\label{propBO}
Let $f\in A^\a (E)$, with $2\a\notin\nn$, and let $\X$ be a
differential
operator of weight $\o$, with $\o<\a$.
Then 
\begeq
|\X(T_x^\a F-T_y^\a F)(z)|\le C\| F\|_\a (d(x,z)+d(y,z))^{\a-\o(\X)}.
 \end{equation}
\end{proposition}

\demo{Proof of theorem \ref{lemabo}}
The proof of theorem 7 is much like the proof of theorem 4, and we
will only sketch it. As a
matter of fact, it is even simpler, as we only have to prove part 2 in
lemma \ref{THELEMMA} to see that $\E(f)$ interpolates $f$, and the
analogous of proposition \ref{smartchange} to see that the function
lies in $A_\a(B)$. Both things can be done at the same time. 

We begin by using a formula for the derivative of the product
analogous to \ref{llit}. Proceeding as there, we see that the terms we
have to bound are of the form:
\begin{equation}
\int_E \X_1\left(\frac1{h_q(z)}\frac{1}{(1-\zbar z)^q}\right)
\X_2\left(\Tz-T^\a_{z_0}F(z)\right)d\mu(\z).
\label{2collons}
\end{equation}
Now we use \ref{cotanucli}, and proposition \ref{propBO}, and from
here on we can proceed as in the proof of theorem 4.\qed

\section{Proof of theorem 8}
To prove the theorem we will need the following:
\begin{proposition}
\label{propbola}
For $\z\in E$, and $a>0$, $R>0$,
$$
\int_{B(\z,R)}d(z,E)^{-a}d\s(z)\le C R^{n-a}
$$
whenever $\Upsilon(E)<n$ and $a<n-\Upsilon(E)$.
\end{proposition}

\demo{Proof}
Let $\eps>0$, and $E_\eps=\{z\in S,\, \rho(z,E)<\eps\}$. 
Just counting the balls of radii $\eps$ contained in $B(\x,R)$ 
implies that for any $n>s>\U(E)$ we have
\begin{equation}
\s(B(\x,R)\cap E_\eps)\le C R^s \eps^{n-s}.    \label{sigma}
\end{equation}
Fix $\x$ and $R$. If $\rho(\x,E)\ge 2\sqrt2 R$, then for any $y\in B(\x,E)$
we have $d(y,E)\ge R$, and the bound is trivial. If 
$\rho(\x,E)\le 2\sqrt2 R$,
then for any $y\in B(\x,R)$ we have that $d(y,E)\le 6R$. Let
$$
B_j=\{z\in B(\x,R),\,\, 6^{-j}R\le \rho(z,E)\le 6^{-j+1} R\}.
$$
We decompose the integral over $B(\x,R)$ into the integrals over
$B_j$, for $j\ge 1$. In each of these integrals we use that
$\rho(x,E)\ge 6^{-j}R$. In this way we see that our integral is
bounded by: 
$$
 \sum_{j=1}^\infty 6^{aj}R^{-a}\s(B_j)\le
\sum_{j=1}^\infty 6^{aj}R^{-a}\s(B(x,R)\cap\{z,\,d(z,E)\le
6^{-j+1}R\}).
$$
If we apply \ref{sigma} to this last sum, and use that $n-s-a>0$ to
see that the sum we obtain is convergent, we are done.   \qed

\demo{Proof of the theorem}
Let $\l=[\be]+1$. We can express $R^\l$ as
$(I+N)^\l=\sum_{k=0}^\l
\binom{\l}{k}N^k$. Now if we write
$$
G_{k,j}(z)=\int_E N^j(\frac1{h_q(z)}\frac1{(1-\zbar z)^q})
N^{k-j}(\T_\z F)(z)\,d\mu(\z)
$$
and
$$
I_{k,j}(z)=\int_0^1
\left(1-t\right)^{\l-\beta-1}
G_{k,j}(tz)\,dt,
$$
we have that
\begin{equation}
\| f\|_{p,1,\be}^p
\le C
\sum_{k=0}^\l \sum_{j=0}^k\int_S |I_{k,j}(rz)|^pd\sigma(z). \label{jotazero}
\end{equation}
The integrals corresponding to $j=0$ can be bounded proceeding in the
same way as we proceeded in the proof of lemma \ref{ojala}. Thus we
will only consider the case $j>0$.

We will write
$$
n_j(z,\z)=N^j\left(\frac1{h_q(z)}\frac1{(1-\zbar z)^q}\right)
$$
so that $0=\int_E n_j(z,\z)\,d\mu(\z)$. Hence, for almost any
$\x\in E$, we have:
$$
G_{k,j}(z)=\int_E n_j(z,\z) N^{k-j}(\T_\z F-\T_\x F)(z)\,d\mu(\z).
$$
Then integrating both sides with respect to 
$(1-\xibar z)^{-q}d\mu(\x)$
and dividing by $h_q(z)$ we get:
$$
G_{k,j}(z)=\frac1{h_q(z)}
\iint_{\ee} \frac{n_j(z,\z)}{(1-\xibar z)^{-q}} 
N^{k-j}(\T_\z F-\T_\x F)(z)\,d\mu(\z)d\mu(\x).
$$
We will write:
$$
I_{k,j}^a(z)=\frac{d(z,E)^{2q-a}}{\mu(B_z)^2}
\iint_{\ee}
\frac{|N^{k-j}(\T_\z F-\T_\x F)(z)|}{d(\z,z)^{q+j-a}d(\x,z)^q}
d\mu(\z)d\mu(\x).
$$
Hence, and applying \ref{cotanucli} and part {\bf (c)} in proposition 
\ref{cotesh}, we have that
$|G_{k,j}(z)|$ is bounded by sums of terms like $I_{k,j}^a(z)$.
On the other hand, $\T_\x F(z)$ is a polynomial in $z$ of degree less
than $2\a$,
so it has a development at $\z$ as
$$
T_\x F(z)=\sum_{|\g|<2\a}\frac1{\g !}D^\g(\T_\x F)(z) w(z,\z)^\g,
$$
wherefrom
\begin{eqnarray*}
\T_\z F(z)-\T_\x F(z)&=&\sum_{\o(\g)<\a}\frac1{\g !}
\left(F_\g(\z)-D^\g (\T_\x F)(\z)\right) w(z,\z)^\g-
\\
&-&\sum_{\a<\o(\g)<2\a}\frac1{\g !} D^\g (\T_\x F)(\z) w(z,\z)^\g.
\end{eqnarray*}
We fix now $k$, $j$, and $a$. Applying to $I_{k,j}^a$ this equality,
 we see that $I_{k,j}^a$ is bounded by sums of terms like
$$
H_{\g,1}(z)=\frac{d(z,E)^{2q-a}}{\mu(B_z)^2}
\iint_{\ee}\frac{|\DELTA||N^{k-j}(w(z,\z)^\g)|}
{d(\z,z)^{q+j-a}d(\x,z)^q}d\mu(\z)d\mu(\x) 
$$
with $\o(\g)<\a$, and
$$
H_{\g,2}(z)=\frac{d(z,E)^{2q-a}}{\mu(B_z)^2}
\iint_{\ee}\frac{|D^\g (\T_\x F)(\z)||N^{k-j}(w(z,\z)^\g)|}
{d(\z,z)^{q+j-a}d(\x,z)^q}d\mu(\z)d\mu(\x) \label{Hgd}
$$
with $\a< \o(\g)<2\a$.

We begin by bounding $H_{\g,2}(z)$. As in lemma \ref{ojala},
$|N^{k-j}(\T_\z f)(z)|$ can be bounded by sums of terms like
$|f_\g(\z)|$.
Using it, \ref{VBB} and part {\bf (a)} in proposition \ref{cotesh}, we
get that
$$
H_{\g,2}(z)\le 
C\frac{d(z,E)^{q+\o(\g)-k}}{\mu(B_z)}\sum_{\o(\de)<\a}
\int_E \frac{|f_\de(\x)|}{d(z,\x)^q}
d\mu(\x).
$$
Recall that we are assuming 
that between $\a$ and $\a+\frac{\u-\la}p$ lies no
integer multiple of $\unmig$, thus $\o(\g)-(\a+\frac{\u-\la}p)>0$. 
Then we can take
$A=\frac{n}p-(\o(\g)-(\a+\frac{s-d}p))+\eps$, with $\eps$ small enough
so that $A$ remains less than $\frac{n}p$.
We write $q=A+q-A$, and apply H\"older's inequality to
this last integral. Of the two integrals we get, we bound the one not
containing $f_\de$ using part {\bf (a)} of proposition \ref{cotesh},
something we can do if $p'(q-A)>\u$, that is, $q>A+\u /p'$. Thus we obtain
$$
\int_E \frac{|F_\de(\x)|}{d(z,\x)^q}d\mu(\x)\le
d(z,E)^{-q+A}\mu(B_z)^{\frac1{p'}}
\left(\int_E
  \frac{|F_\de(\x)|^p}{d(z,\x)^{Ap}}d\mu(\x)\right)^{\frac1p}.
$$
Then, if we use $U_\u'$ to bound $\mu(B_z)$, and substitute $k$ by $\l$
(something we can do as $k\le \l$), we get that 
$$
H_{\g,2}(z) \le
C d(z,E)^{A+\o(\g)-\l-\frac{\u}{p}}\sum_{\o(\de)<\a}
\left(\int_E
  \frac{|F_\de(\x)|^p}{d(z,\x)^{Ap}}d\mu(\x)\right)^{\frac1p}.
$$
If we apply this to bound $\int_0^1 (1-t)^{\l-\be-1}H_{\g,2}(tz)dt$ we
get that it is bounded by
$$
 C\sum_{\o(\de)<\a}\int_0^1
(1-t)^{\l-\be-1}
d(tz,E)^{A+\o(\g)-\l-\frac{\u}{p}}dt
\left(\int_E\frac{|F_\de(\x)|^p}{d(z,\x)^{Ap}}d\mu(\x)\right)^{\frac1p}.
$$
 Now in the integral with respect to $t$, the sum of the exponents
 remains positive, by the way we chose $A$. Thus it is bounded by some
constant depending on the exponents, but not on $z$. 
If we use this bound when integrating over $S$, we get that the
 integral of $\Hgd$ is bounded by  
$$
 \sum_{\o(\de)<\a}\int_S
\int_E\frac{|f_\de(\x)|^p}{d(z,\x)^{Ap}}d\mu(\x)d\s(z)\le
\sum_{\o(\de)<\a}\|f_{\de}\||{\Lpmu},
$$
as $Ap<n$.

\smallskip

Now we have to bound $\Hgu(z)$.
We use again \ref{VBB}, and  that $k\le \l$, and get
$$
\Hgu(z)\le
\frac{d(z,E)^{2q-a}}{\mu(B_z)^2}
\iint_{\ee}\frac{|\DELTA|}
{d(\z,z)^{q+\l-a-\o(\g)}d(\x,z)^q}d\mu(\z)d\mu(\x).
$$
As a matter of fact, we have to bound
$(\int_0^1(1-t)^{\l-\be-1}\Hgu(tz)dt)^p$,
which is thus bounded by:
$$
\left(
\int_0^1(1-t)^{\l-\be-1}
\frac{d(tz,E)^{2q-a}}{\mu(B_{tz})^2}
\iint_{\ee}\frac{|\DELTA|d\mu(\z)d\mu(\x)}
{d(\z,tz)^{q+\l-a-\o(\g)}d(\x,tz)^q}dt
\right)^p.
$$
Let $0<\de<\min\{\l-\be,n-\U(E)\}$be small enough. 
Let $A=\a-\o(\g)+(2\u-\la)/p+\de$. 
Then because of H\"older's inequality, the previous integral is
bounded by:
$$
\int_0^1(1-t)^{(\l-\be-\de)p-1}
\frac{d(tz,E)^{(A+3\de) p}}{\mu(B_{tz})^2}
\iint_{\ee}\frac{|\DELTA|^p d\mu(\z)d\mu(\x)\,dt}
{d(\z,tz)^{(\l-\o(\g)+\de)p}d(\x,tz)^{Ap}}
$$
multiplied by
\begin{equation}
\left(\int_0^1\iint_{\ee}
\frac{(1-t)^{\de p'-1} d(tz,E)^{(2q-a-A-3\de)p'}}
{\mu(B_{tz})^2 d(\z,tz)^{(q-a-\de)p'}d(\x,tz)^{(q-A)p'}}
d\mu(\z)d\mu(\x)\,dt
\right)^{\frac{p}{p'}}. \label{uau}
\end{equation}
Assume that $q>a+\de+\u/p'$ and that $q>A+\u /p'$. Then we can apply
part {\bf (a)} of proposition \ref{cotesh} with the inner integral
in \ref{uau}. If we use lemma \ref{lemashiti} to the integral
with respect to $t$ we obtain, we see that \ref{uau} is bounded by 
 $d(z,E)^{-\de p'}$, so that
$\int_S \left(
\int_0^1(1-t)^{\l-\be-1} \Hgu(tz)\,dt
\right)^p d\s(z)$
is, because of Fubinni's theorem, bounded by
$$
\iint_{\ee}|\DELTA|^p J(\x,\z)
d\mu(\z)d\mu(\x)
$$
where
$$
J(\x,\z)=
\int_S d(z,E)^{-\de p}\int_0^1
\frac{(1-t)^{(\l-\be-\de)p-1} d(tz,E)^{(A+3\de)p}}
{\mu(B_{tz})^2 d(\z,tz)^{(\l-\o(\g)+\de)p}d(\x,tz)^{A p}}dt\,d\s(z).
$$
Thus to finish the proof we only have to show that $J(\x,\z)$ can be
bounded by $d(\z,\x)^{-(\a-\o(\g))+\la}\mu[\x,\z]^{-2}$.

To bound $J(\x,\z)$ we split $S$ into the sets 
$S_1(\z,\x)=\{ z\in S,\,d(z,\z)\le d(z,\x)\}$, and $S_2=S\setminus
S_1$. We will consider only the integral over $S_1$ as that over $S_2$
is done likewise. 

Let $(tz)_0$ be a point in $E$ nearest to $tz$. Then for 
$x\in B(\x,d(\x,tz))$, we have 
$d(x,(tz)_0)\le 6 d(\x,tz)$. On $S_1$ we also have that
$d(\z,\x)\le 4 d(tz,\x)$. Using these facts and proceeding as 
for \ref{cotamuBx}, we get that
$$
\mu(B_{tz})\ge C\frac{\rho(tz,E)^\u}
{\rho(\x,\z)^\la \rho(\x,tz)^{\u-\la}}\mu[\x,\z].
$$
Thus, the part of $J(\x,\z)$ corresponding to the integral over $S_1$
is bounded by:
$$
\frac{d(\z,\x)^{2\la}}{\mu[\z,\x]^2}
\int_S d(z,E)^{-\de p}\int_0^1
\frac{(1-t)^{(\l-\be-\de)p-1} d(tz,E)^{(A+3\de)p-2\u}}
{ d(\z,tz)^{(\l-\o(\g)+\de)p}d(\x,tz)^{A p-2\u+2\la}}dt\,d\s(z).
$$
To bound this integral, we split $S_1$ into $B_1=S_1\cap
B(\x,d(\z,\x)/4)$ and its complementary. 

On $B_1$, we use that $d(\z,\x)\le 4d(tz,\x)$ and that $d(tz,E)\le
d(tz,\z)$, together with the value of $A$, to bound the integral by
$$
\frac{d(\z,\x)^{-Ap+2\u}}{\mu[\z,\x]^2}
\int_S d(z,E)^{-\de p}
\int_0^1
\frac{(1-t)^{(\l-\be-\de)p-1} }
{ d(\z,tz)^{(\l-\a-3\de)p+\la} }dt\,d\s(z).
$$
Now applying lemma \ref{lemashiti} to the inner  integral we get that
the previous integral is bounded by:
$$
\frac{d(\z,\x)^{-Ap+2\u}}{\mu[\z,\x]^2}
\int_S 
\frac{ d(z,E)^{-\de p} }
{ d(\z,z)^{n-2\de p} }d\s(z).
$$
Now we split $B_1$ into the sets $B^j=B(\z,2^{-j}d(\z,\x)/4)$ for
$j\ge 1$. Then
we decompose the integral over $B_1$ into the integrals over
$B^{j+1}\setminus B^{j}$. In each of these integrals de use that
$d(z,\z)\ge 2^{j-2}d(\z,\x)$ and apply proposition \ref{propbola} to
each of them, we get that our integral is bounded by $d(\z,E)^{\de p}$
multiplied by a finite sum,
and so we are done.

To bound the integral over $S_1\setminus B_1$, we use that
$d(tz,\x)\ge d(z,\z)$, and thus what we have to bound is:
$$
\frac{d(\z,\x)^{2\la}}{\mu[\z,\x]^2}
\int_S d(z,E)^{-\de p}\int_0^1
\frac{(1-t)^{(\l-\be-\de)p-1} d(tz,E)^{(A+3\de)p-2\u}}
{ d(\z,tz)^{(\l-\o(\g)+\de)p}d(\z,z)^{A p-2\u+2\la}}dt\,d\s(z).
$$
Now proceeding as before, but using the descomposition of
$S_1\setminus B_1$ into $B^j$ with $j<0$, we obtain again the same bound.

\smallskip

It remains to see that $\E(f)$ interpolates $f$. It is shown in
essentially the same way as in lemma \ref{viscaelbarsa},
 we will only sketch the differences. 

We first observe that $D^\g f(z)$ can be written as
$$
D^\g \E(f)(z)=\sum_{\o(\X')+\o(\X'')=\o(\X)} \int_E \X'(n(z,y))\,
\X''(\T_y f) (z)\, d\mu(y).
$$
In this sum, the only term with $\o(\X')=0$ has $\X''=D^\g$ and
$\X'=Id$. All the other terms have $\o(\X')>0$, whence 
$\int_E \X'(n(z,\cdot))d\mu=0$. This allows us to add to these terms 
anything not depending on $y$. Thus the difference 
$D^\g \E(f)(z)-f_\g (\z)$ can be expressed as 
\begin{equation}
\int_E n_0(z,y) \left[ D^\g \T_y f(z)- f_\g (\z)\right] d\mu(y)   \label{head}
\end{equation}
plus terms of the form 
\begin{equation}
\int_E \X'(n_0(z,y)) \X''( \T_y f(z)- \T_x f (z)) d\mu(y)   \label{tails}
\end{equation}
for any $x\in E$ for which $\T_x f$ is well defined. 

To bound the terms in \ref{tails} we develop $\T_x f (\cdot)$ at the
point $y$. As it is a polynomial in $z$ of degree $\le 2\a$, two kinds
of terms arise: those with weight less than $\a$ and those with
weight greater than $\a$. The terms with weight less than $\a$ are
handled in the same way as in lemma \ref{viscaelbarsa}.

The terms with weight greater than $\a$ have the form 
$D^a (\T_x f)(y) \X''(w(y,z)^a)$, with $\o(a)>\a$. Using
\ref{cotanucli} and that $|\X(w(y,z)^a)|\le d(y,z)^{\o(a)-\o(\X)}$, we
see that these are bounded by sums of $d(z,E)^{\o(a)-\o(\X)}|f_b(x)|$,
with $\o(b)<\a$. Then we integrate it with respect to
$d(z,x)^{-q}d\mu(x)$ and divide it by $h_q(z)$ (we must do it, as we
did the same for the terms with weight less than $\a$), and finish as
in the proof of lemma \ref{ojala}.

It remains to bound \ref{head}. We first observe that 
$f_\g(\z)=D^\g (\T_\z f)(\z)$. Moreover, for those $\z$ for which 
$\T_\z f$ is well defined, $D^\g (\T_\z f)(\cdot)$ is a polynomial,
so $|D^\g (\T_\z f)(z)-D^\g (\T_\z f)(\z)|$ is bounded by 
$d(z,\z) h(\z)$, for some $h\in L^p (\mu)$. Hence it is enough to
bound
$$
\int_E n_0(z,y)  D^\g (\T_\z f- \T_y f) (\z) d\mu(y).
$$
As before, we develop $\T_\z f (\cdot)$ at the point $y$. Then two
kinds of terms arise: the first terms can be reagrupated in
$$
\int_E n_0(z,y) \D_a (\z,y) D^\g(w(y,z)^a) d\mu(y),
$$
with $\o(a)<\a$, and then spare terms like
$$
\int_E n(z,y) \D_a (\z,y) D^\g(w(y,z)^a) d\mu(y),
$$
with $\o(a)>\a$. The former term is dealt with as in the proof of
lemma \ref{viscaelbarsa}. The latter terms are easily bounded by sums
of $|f_b(\z)| d(z,\z)^{\o(\a)-\o(\X)}$, thereby approaching $0$ when
$d(z,\z)\to 0$.
\qed

\end{document}